\documentclass[preprint,12pt]{amsart}

\usepackage{amssymb}
\usepackage{mathrsfs}   

\usepackage{tikz}

\usepackage[hang,small,bf]{caption}
\newcommand\mL{L\kern-0.08cm\char39}  

\newtheorem{theorem}{Theorem}[section]
\newtheorem{corollary}[theorem]{Corollary}
\newtheorem{lemma}[theorem]{Lemma} 
\newtheorem{definition}[theorem]{Definition}

\theoremstyle{definition}
\newtheorem{exmp}[theorem]{Example}  

\def\proof{\noindent {\bf Proof. }}
\def\cbd{\hfill$\Box$}   

\def\f{{\it \overline f}}
\def\g{{\it \overline g}}

\begin{document}

\title[On generic and dense chaos on hyperspaces]{On generic and dense chaos for maps induced on hyperspaces}



\author[M. Ml\'{\i}chov\'a]{Michaela Ml\'{\i}chov\'a}
\email{Michaela.Mlichova@math.slu.cz}

\author[M.  \v Stef\' ankov\' a]{Marta \v Stef\' ankov\' a}
\email{Marta.Stefankova@math.slu.cz}

\address{Mathematical Institute, Silesian University, 746 01 Opava, Czech Republic}

\begin{abstract}
A continuous map $f$ on a compact metric space $X$ induces in a natural way the map $\f$ on the hyperspace 
$\mathcal K(X)$  of all closed non-empty subsets of $X$. We study the question of transmission of chaos between $f$ and $\f$. We deal with generic, generic $\varepsilon$-, dense and dense $\varepsilon$-chaos for interval maps. We prove that 
all four types of chaos transmit from $f$ to $\f$, while the converse transmission from $\f$ to $f$ is true for
generic, generic $\varepsilon$- and dense $\varepsilon$-chaos. Moreover, the transmission of dense $\varepsilon$-
and generic $\varepsilon$-chaos from $\f$ to $f$ is true for maps on general compact metric spaces.\\
\end{abstract}

\maketitle
\section{Introduction and preliminary results}

Let $(X,d)$ be a compact metric space endowed with the metric $d$.
A dynamical system $(X,f)$, where $f: X \to X$ is continuous,
induces in a natural way the system $(\mathcal K(X), \f)$ on the hyperspace $\mathcal K(X)$ consisting of all closed non-empty subsets of $X$.
A natural question arises what are the connections between the {\em individual} dynamics given by $(X,f)$ and the {\em collective} dynamics given by $(\mathcal K(X), \f)$.

\subsection{Hyperspaces}

Let us recall some definitions and properties concerning hyperspaces. We will give here only those we will need in the following text. For the proofs and further study of topology of hyperspaces see, e.g., \cite{EM}, \cite{IN} and \cite{Mac}. 

The map $\f : \mathcal K(X) \to \mathcal K(X)$ induced by a continuous map $f :X \to X$ on the space 
$\mathcal K(X) = \{ K\subset X;\ K \neq\emptyset  {\rm \ is\ compact} \}$ is defined by
$\f(K) = f(K) = \{ f(x);\ x\in K \}$, $K\in\mathcal K(X)$; note that $\mathcal K(X)$ is invariant for $\f$,
and if $f$ is continuous then $\f$ is also continuous.

Let  $x\in X$ be a point, $A\subset X$ a non-empty set, and $\varepsilon > 0$. Define the {\em distance from the point $x$ to the set $A$} by
${\rm dist}(x, A) := {\rm dist}(\{ x\}, A) = \inf\{d(x, y);\  y\in A \}$,
and the {\em $\varepsilon$-neighborhood of the set $A$} by
$N(A, \varepsilon) := \{ x\in X;\  {\rm dist}(x, A)<\varepsilon\}$.
The {\em Hausdorff distance} on  $\mathcal K(X)$  is defined as follows.
For any $A, B \in \mathcal K(X)$,
$$
d_H(A, B) := \inf \{ \varepsilon \ge 0;\  A\subset N(B, \varepsilon) {\rm \ and\ } B\subset N(A,\varepsilon) \}.
$$
Let us note that $(\mathcal K(X), d_H)$ is a compact metric space.

Denote by $\mathbb N$  the set of positive integers.
Let $S_1,\dots , S_n \subset X$, $n\in \mathbb N$, be a finite collection of non-empty sets. Define a subset
$\langle S_1, \dots , S_n \rangle$ of $\mathcal K(X)$ by
$$
\langle S_1, \dots , S_n \rangle := \{ K\in\mathcal K(X);\  K\subset\bigcup_{i=1}^n S_i,\ K\cap S_i \neq\emptyset {\rm \ for\ }
i=1, \dots , n \}.
$$

\noindent
The family 
$$
\mathcal B = \{\langle U_1, \dots , U_n \rangle ;\ U_1,\dots , U_n {\rm \ are\ open\ in\ } X {\rm\ and\ } n\in\mathbb N \} 
$$

\noindent
of subsets of $\mathcal K(X)$ is a basis for the so called {\em Vietoris topology} on $\mathcal K(X)$. Note that the topology induced by the Hausdorff metric and the Vietoris topology for $\mathcal K(X)$ coincide.
%


The first systematic investigation of dynamical properties of induced maps was done by Bauer and Sigmund in their paper 
``Topological dynamics of transformations induced on the space of probability measures'' from 1975 (see \cite{BauSig}). 
The authors studied which of topological properties of the individual system $(X,f)$ (like, e.g., transitivity, mixing properties, specification properties, distality) are ``inherited'' by collective system $(\mathcal K(X), \f)$. Note that, as the title suggests,  the authors in this paper considered simultaneously systems induced on the space of probability measures.

There are  also many papers concerning connections between chaotic behavior of individual and collective systems. Properties like Devaney chaos, Li-Yorke chaos, distributional chaos, $\omega$-chaos or topological chaos have been considered (see, e.g., \cite{GKLOP} where also an extensive list of references is given).
The main aim of the present paper is to study connections between other two kinds of chaos, namely generic and dense chaos, of $f$ and $\f$.

\subsection{Generic and dense chaos}\label{chaos}

Let $f\in C(X)$, the class of continuous maps $X\to X$,  and let $\varepsilon > 0$.
In the following, we will use the notation which is taken from \cite{S1}:
$$C_1(f) := \{ (x, y) \in X^2;\ \liminf_{n\to\infty} d(f^n(x),  f^n(y)) = 0 \},$$
$$C_2(f) := \{ (x, y) \in X^2;\ \limsup_{n\to\infty} d(f^n(x),  f^n(y)) > 0 \},$$
$$C_2(f, \varepsilon) := \{ (x, y) \in X^2;\ \limsup_{n\to\infty} d(f^n(x),  f^n(y)) > \varepsilon \},$$
$$C(f) := C_1(f) \cap  C_2(f), $$
$$C(f, \varepsilon) := C_1(f) \cap C_2(f, \varepsilon).$$

\noindent

Note that $C(f)$ (resp. $C(f, \varepsilon)$) is the set of {\em Li-Yorke pairs} (resp. {\em $\varepsilon$-Li-Yorke pairs}.
In this notation, the well known definition of chaos in the sense of Li and Yorke reads as follows: A map $f\in C(X)$ is {\em LY-chaotic} if there exists an uncountable set $S$ such that
$C(f) \supset S\times S \setminus \{ (x, x);\ x\in X \}$.

In the 80s, A. Lasota proposed a new concept of Li-Yorke chaos, the so called generic chaos, see \cite{Pio1}.
Inspired by this idea, \mL . Snoha \cite{S1} introduced three other variants of this notion. Recall that a set is 
of {\em first category} if it is a union of a countable family of nowhere dense sets; a set that is not of first category is a {\em second category} set.
A set is {\em residual} if its complement is a first category set or, equivalently, if it contains a dense $G_{\delta}$ subset.

\begin{definition}\label{def}
A map $f\in C(X)$ is called
\begin{enumerate}
\item [(i)]
	{\em generically chaotic}, if the set $C(f)$ is residual in $X^2$,
\item [(ii)]
	{\em generically $\varepsilon$-chaotic}, if the set $C(f,\varepsilon)$ is residual in $X^2$,
\item [(iii)]
	{\em densely chaotic}, if the set $C(f)$ is dense in $X^2$,
\item [(iv)]
	{\em densely $\varepsilon$-chaotic}, if the set $C(f,\varepsilon)$ is dense in $X^2$.
\end{enumerate}
\end{definition}

Properties of generic chaos for interval maps have been very deeply studied in \cite{S1}. 
Snoha showed there that generically chaotic, generically $\varepsilon$-chaotic and densely $\varepsilon$-chaotic interval maps, where $\varepsilon$ is a positive real number, are equivalent. Moreover, he characterized such maps in terms of behavior of subintervals of $I = [0,1]$. Densely chaotic interval maps have been characterized in  \cite{S2}. It has been proved there, among others, that in the class of piecewise monotone maps with finite number of pieces, dense chaos and generic chaos coincide. 
Properties of generically and densely chaotic interval maps are summarized in Theorems \ref{th1} and \ref{th2}.
Note that the minimum of topological entropies of both generically and densely chaotic interval maps is $(1/2)\log 2$, see \cite{S1} and \cite{Rue}.

Properties of generically $\varepsilon$-chaotic maps on metric spaces have been studied in \cite{Mur}. Murinov\' a showed there that many of the properties of generically $\varepsilon$-chaotic interval maps from \cite{S1} can be extended to a large class of metric spaces, e.g., that generic $\varepsilon$-chaos is equivalent to dense $\varepsilon$-chaos
 but, on the other hand, there is a convex continuum on the plane on which generic chaos and generic $\varepsilon$-chaos are not equivalent.

Let us recall here characterizations of generically (resp. densely) chaotic interval maps proved by Snoha.
Note that Theorem \ref{th1} is rewritten from \cite{S1}, since we use almost all properties occurring there (the only properties we do not need are (d) and (e) but they are interesting in themselves). The properties in Theorem \ref{th2} are chosen from \cite{S2}, Theorem 1.2, and the subsequent text in the way most convenient for our purposes.

 Let $I$ be the compact unit interval $[0, 1]$.
 By an interval we mean a nondegenerate (not necessarily compact) interval lying in $I$. If $J$ is an interval then diam\,$J$ denotes its length. The {\it distance of two sets} $A, B \subset I$ is defined by dist$(A,B) := \inf \{ |x-y|;\ x\in A,\ y\in B \}$, and recall that the distance from a point $a$ to a set $B$ is dist$(a,B) =\ $dist$(\{ a\},B)$.
A compact interval $J$ is an {\em invariant transitive interval of} $f$ if  $f(J)\subset J$ and the restriction of $f$ to $J$ is topologically transitive. For any set $A$, int$(A)$ (resp. cl$(A)$) denotes the interior (resp. closure) of $A$.
Denote Orb$(f, A) := \bigcup_{n=0}^{\infty} f^n (A)$.

\begin{theorem} \label{th1} {\rm({\cite{S1},  Theorem 1.2)}}
Let $f \in C(I)$. The following conditions are equivalent:

\begin{enumerate}

\item [(a)] $f$ is generically chaotic,

\item [(b)] for some $\varepsilon>0$, $f$ is generically $\varepsilon$-chaotic,

\item [(c)] for some $\varepsilon>0$, $f$ is densely $\varepsilon$-chaotic,

\item [(d)] $C_1(f)$ is dense in $I^2$ and $C_2(f)$ is a second category set in any interval $J^2 \subset I^2$,

\item [(e)] $C_1(f)$ is dense in $I^2$ and for some $\varepsilon >0$, $C_2(f, \varepsilon)$ is dense in $I^2$, 

\item [(f)] the following two conditions are fulfilled simultaneously:
	\begin{enumerate}
	
	\item [(f-1)] for every two intervals $J_1, J_2$, $\liminf \limits_{n \to \infty} {\rm dist} (f^n(J_1), f^n(J_2))=0$,
	
	\item [(f-2)] there exists an $a>0$ such that for every interval $J$, \\
	$\limsup \limits _{n \to \infty} {\rm diam\,} f^n (J)>a $,
	
	\end{enumerate}

\item [(g)] the following two conditions are fulfilled simultaneously:
	\begin{enumerate}
	
	\item [(g-1)]  there exists a fixed point $x_0$ of $f$ such that for every interval $J$, $\lim \limits_{n \to \infty} {\rm dist} (f^n(J), x_0)=0$,
	
	\item [(g-2)] there exists a $b>0$ such that for every interval $J$, \\
	$\liminf \limits _{n \to \infty} {\rm diam\,} f^n (J)>b$,
	
	\end{enumerate} 

\item [(h)] the following two conditions are fulfilled simultaneously:
	\begin{enumerate}
	
	\item [(h-1)] $f$ has a unique invariant transitive interval or two invariant transitive intervals having one point in common, 
	
	\item [(h-2)]  for every interval $J$ there is an invariant transitive interval $T$ of $f$ such that ${\rm Orb} (f, J) \cap {\rm int} (T) \neq \emptyset$.
	
	\end{enumerate}

\end{enumerate}
Moreover, the equivalences ${\rm (b)} \Leftrightarrow {\rm (c)}  \Leftrightarrow {\rm (e)}  \Leftrightarrow {\rm (f)} $ hold with the same $\varepsilon$ and with $a = \varepsilon$ in {\rm (f-2)}.

\end{theorem}

 \begin{theorem}\label{th2} \cite{S2}
 A function $f \in C(I)$ is densely chaotic if and only if the following three conditions are fulfilled simultaneously:
 \begin{enumerate}
 
 \item[(a)] there is a fixed point $x_0$ of $f$ such that for every interval $J$, 
 $$\lim\limits_{n \to \infty} {\rm dist} (f^n(J), x_0)=0,$$
 
 \item[(b)] for every interval $J$, $\liminf_{n\to \infty} {\rm diam\,} f^n (J)>0$,

 \item[(c)] every one-sided punctured neighbourhood of the point $x_0$ contains points $x, y$ with $(x, y) \in C(f)$ and moreover, if $x_0 \in {\rm int} (I)$ then every neighbourhood of $x_0$ contains points $x<x_0<y$ with $(x, y)\in C(f)$.
 
 \end{enumerate}
 \end{theorem}

In Section 2 we show that, for interval maps, if $\f$ is generically chaotic then also $f$ is generically chaotic. In Section 3 we prove that, for maps on general compact metric spaces, dense $\varepsilon$-chaoticity of $\f$ implies that $f$ has this property, too.  Section 4 concerns the opposite implications; we show that, for interval maps, dense (resp., dense 
$\varepsilon$-) chaos of $f$ imply dense (resp., dense  $\varepsilon$-) chaos of $\f$. In Section 5 we give some results concerning the question of transmission of dense chaos from $\f$ to $f$. We also provide a scheme, where the obtained results together with their corollaries are given in a clearly arranged form.


\section{Generically chaotic $\f$ implies generically chaotic $f$}

In this section we will prove the following
\begin{theorem}\label{main1}
Let $f \in C(I)$ be such that the induced map $\f$ is generically chaotic. Then the function $f$ is also  generically chaotic.
\end{theorem}

 
%
%
%
%
%
%

To prove Theorem \ref{main1} we need several lemmas. The next lemma, which follows easily from uniform continuity of $\f$, is a version of Lemmas 4.1 and 4.2 from \cite{S1}. 

\begin{lemma}\label{new1}
Let $f \in C(I)$, $g=f^k$ for some positive integer $k$, and let $\f$ and $\g$ be induced by $f$ and $g$. \noindent Moreover, let $\mathcal{A, B} \subset \mathcal{K}(I)$ be non-empty sets. Then

\begin{enumerate}

\item[(i)] 
	$\liminf\limits_{ n \to \infty}{\rm dist} (\f^n(\mathcal A), \f^n(\mathcal B))=0$ iff $\liminf\limits_{ n \to \infty}{\rm dist} (\g^n(\mathcal A), \g^n(\mathcal B))=~0$,

\item[(ii)] 
	$\lim\limits_{ n \to \infty}{\rm dist} (\f^n(\mathcal A), \f^n(\mathcal B))=0$  iff $\lim\limits_{ n \to \infty}{\rm dist} (\g^n(\mathcal A), \g^n(\mathcal B))=~0$,  

\item[(iii)] 
	$\liminf \limits _{n \to \infty} {\rm diam\,} \f^n (\mathcal A)=0$ iff $\liminf \limits _{n \to \infty} {\rm diam\,} \g^n (\mathcal A)=0$,

\item[(iv)]  
	$\limsup \limits _{n \to \infty} {\rm diam\,} \f^n (\mathcal A)=0$ iff $\limsup \limits _{n \to \infty} {\rm diam\,} \g^n (\mathcal A)=0$,
\item [(v)] 
	$C_1(\f)=C_1(\g)$ and $C_2(\f)=C_2(\g)$,

\item [(vi)] 
	$\f$ is generically or densely chaotic if and only if $\g$ is generically or densely chaotic.

\end{enumerate}
\end{lemma}

%
%

\begin{lemma} \label{lm3}{\rm (\cite{S1}, Lemma 4.3)}
Let $f \in C(I)$. Then the following three conditions are equivalent:

\begin{enumerate}
\item [(i)] 
	$C_1(f)$ is residual in $I \times I$,

\item [(ii)] 
	$C_1(f)$ is dense in $I \times I$,

\item [(iii)] 
	for every two intervals $J_1, J_2$, $\liminf \limits_{n \to \infty} {\rm dist} (f^n(J_1), f^n(J_2))=0$ (i.e., condition (f-1) from Theorem \ref{th1}).
\end{enumerate}
\end{lemma}

\begin{lemma}\label{new3}
Let $f \in C(I)$ and let $\f$ be induced by $f$. If $C_1(\f)$ is residual in $\mathcal K(I) \times \mathcal K(I)$, then $C_1(f)$ is residual in $I \times I$.
\end{lemma}

\proof Let $J_1$ and $J_2$ be arbitrary intervals. Since $C_1(\f)$ is residual, there exist non-empty closed sets $U_1 \subset J_1$, $U_2 \subset J_2$ such that 
$$\liminf\limits_{n \to \infty} d_H (\f ^n (U_1), \f ^n(U_2))=0.$$ Obviously, $\liminf\limits_{n\to \infty}{\rm dist}(f ^n (U_1), f ^n(U_2))=0$ and thus $f$ satisfies {\rm (iii)}  from Lemma \ref{lm3}. 
\cbd \\

The following property is easy.

\begin{lemma}\label{top}
In arbitrary topological space, if $B \subset A$ is dense in $A$, then $\mathcal K(B)$ is dense in $\mathcal K(A)$.
\end{lemma}

The first part of the proof of the following lemma might seem to be identical to the proof of Lemma 4.8 in \cite{S1}. 
But, in fact, our lemma has different assumptions,  in the proof we sometimes work with $\f$ (not with the original $f$)
and, moreover, in the conclusion of our proof we need some sets constructed in the first part. Hence
we give here the full version of this proof.

\begin{lemma} \label{lm8} 
Let $f \in C(I)$ and let $\f$ be induced by $f$. Let $x_0 \in I$ be a fixed point of $f$ and let $\f$ be generically chaotic. Then there exists a $\delta>0$ such that no interval containing $x_0$ and with diameter less than $\delta$ is $f$-invariant.
\end{lemma}

\proof 
Since the closure of an invariant interval is an invariant interval with the same diameter, it suffices to prove the claim of our lemma for compact intervals. Assume on the contrary that for every $\delta >0$ there is a compact invariant interval $J(\delta)$ containing $x_0$ and with diameter less than $\delta$. Then infinitely many of the intervals $J(1/n)$, $n =1, 2, \ldots$, have the right endpoints greater than $x_0$ or infinitely many of them have the left endpoints less than $x_0$. Without loss of generality we may suppose the first possibility. Further, observe that the intersection of two compact invariant intervals is a compact invariant interval. Now it is not difficult to see that there exists a sequence of invariant intervals $J_n=[x_0-a_n, x_0+b_n]$, $n=1, 2, \ldots$, where $\lim_{n\to \infty}a_n=0$, $\lim_{n\to \infty}b_n=0$, and for every $n$, $0<b_{n+1}<b_n, 0 \le a_{n+1} \le a_n$, and $a_{n+1}=a_n$ if and only if $a_n=0$. Consider two cases.

\noindent {\bf Case 1.} For every $n$, $a_n>0$. Then for every $n$ we have $J_{n+1} \subset {\rm int} J_n$. Let $m$ be a positive integer. By Lemmas \ref{new3} and \ref{lm3}, the set $C_1(f)$ is residual in $I \times J_{m+1}$. Thus there exists a set $B_m \subset I$ such that $B_m$ is residual in $I$ and for every $x \in B_m$ there is $y \in J_{m+1}$ with $\liminf_{n \to \infty} |f^n(x)-f^n(y)|=0$. Since $J_{m+1} \subset {\rm int} J_m$ and the interval $J_m$ is invariant, we can see that for every $x \in B_m$ there exists a positive integer $m(x)$ such that ${\rm Orb}(f^{m(x)}(x))\subset J_m$. Now we consider the set $B=\bigcap_{m=1}^\infty B_m$. It is residual in $I$ and it is easy to see that for every $x \in B$, $\lim_{n \to \infty} f^n(x)=x_0$.

\noindent {\bf Case 2.} For some $n$, $a_n=0$. Without loss of generality we may assume that $a_1=0$ and consequently $a_n=0$ for all $n$. Now we cannot use the inclusions $J_{n+1} \subset {\rm int} J_n$ from Case 1. But it suffices to take into account that $f(J_1) \subset J_1$, and analogously as in Case 1 we can provide that the orbits of points from $J_1$ generically converge to $x_0$, i.e., that there exists a set $B \subset J_1$, $B$ residual in $J_1$, such that $\lim_{n \to \infty} f^n(x)=x_0$ for all $x\in B$.

We can see that in either case there are an interval $A \subset I$ and a set $B \subset A$, $B$ is residual in $A$, such that for any $x \in B$, $\lim_{n \to \infty}f^n(x)=x_0$. So $B$ contains an intersection of countably many open dense sets $G_n$, and hence $\mathcal K(B) \supset \mathcal K\left(\bigcap_{n=1}^\infty G_n\right)$. By the definition of $\mathcal K(\cdot)$ we have

\begin{eqnarray*}
\mathcal K\left(\bigcap_{n=1}^\infty G_n\right)  &= &\{P \in \mathcal K(A): P \subset G_n\ {\rm for\ any}\ n\} = \\
& = & \bigcap_{n=1}^\infty \{P \in \mathcal K(A): P \subset G_n\}= \bigcap_{n=1}^\infty \langle G_n\rangle.
\end{eqnarray*}

Since the sets $\langle G_n\rangle$ are open and, by Lemma \ref{top}, they are dense in $\mathcal K(A)$, the set $\mathcal K (B )$ is residual in $\mathcal K(A)$. Obviously, $\mathcal K(B) \times \mathcal K(B)$ is residual in $\mathcal K(A) \times \mathcal K(A)$, and thus $C_2(\f)$ is of first category in $\mathcal K(A) \times \mathcal K(A)$. This contradicts our assumption that $\f$ is generically chaotic (i.e., $C_2(\f)$ is residual).
\cbd
\bigskip

The following three lemmas from \cite{S1} will be used in the proof of the subsequent Lemma \ref{lm16}.

\begin{lemma}\label{lm7}{\rm (\cite{S1}, Lemma 4.7)}
Let $f\in C(I)$ and $J$ be a compact interval with $\limsup_{n\to \infty} {\rm diam\,} f^n(J)>0$. Then ${\rm Orb}(f, J)$ contains a periodic point of $f$. Moreover, if the conditions (f-1) from Theorem \ref{th1} is fulfilled then ${\rm Orb}(f, J)$ contains a periodic point of $f$ with period $1$ or $2$.
\end{lemma}

\begin{lemma}\label{lm9}{\rm(\cite{S1}, Lemma 4.9)}
Let $f\in C(I)$, $x_0$ be a fixed point of $f$ and let  for every two intervals $J_1, J_2$, $\liminf \limits_{n \to \infty} {\rm dist} (f^n(J_1), f^n(J_2))=0$ (i.e., the condition (f-1) from Theorem \ref{th1} be fulfilled). Let there exist arbitrarily small $f$-invariant intervals arbitrarily close to the point $x_0$. Then there exist arbitrarily small $f$-invariant intervals containing the point $x_0$.  
\end{lemma}

\begin{lemma}\label{lm12}{\rm(\cite{S1}, Lemma 4.12)} Let $f \in C(I)$ and $J$ be an interval. Then the following two conditions are equivalent:

\begin{enumerate}
\item[(i)] 
	there are $x, y \in J$ with $\liminf\limits_{n \to \infty} |f^n(x)-f^n(y)|>0$,

\item[(ii)] 
	there are $x, y \in J$ with $\limsup\limits_{n \to \infty} |f^n(x)-f^n(y)|>0$.
\end{enumerate} 

\noindent Further, the following two conditions are equivalent:

\begin{enumerate}
\item[(iii)]  
	$\liminf\limits_{n \to \infty} {\rm diam\,} f^n(J)>0$,

\item[(iv)] 
	$\limsup\limits_{n \to \infty} {\rm diam\,} f^n(J)>0$,
\end{enumerate} 

\noindent and either of the conditions {\rm (i)}, {\rm (ii)} implies either of the conditions {\rm (iii)},\ {\rm (iv)}. Moreover, if $J$ is a compact interval then all the conditions {\rm (i)} -- {\rm (iv)} are equivalent.
\end{lemma}

The following lemma together with its proof is based on \cite{S1}, Lemma 4.16, the part concerning implication (vi) \& (f-1) $\Rightarrow$ (i).

\begin{lemma}\label{lm16}
Let $f \in C(I)$ and let $\f$ be induced by $f$. Assume that $\f$ is generically chaotic. Then there exists a real number $a>0$ such that for any interval $J \subset I$, $\limsup_{n \to \infty} {\rm diam\, } f^n(J)>a$, i.e., $f$ satisfies the condition (f-2) from Theorem \ref{th1}.
\end{lemma}

\proof By \cite{S1}, Lemma 4.1 (an original version of our Lemma \ref{new1} for $f$), it suffices to prove the claim of our lemma for $g=f^2$,  i.e., there exists $a>0$ such that for every interval $J \subset I$, $\limsup_{n \to \infty} {\rm d iam} g^n(J)>a$.  Moreover, we can consider only the intervals containing fixed points of $g$.  Indeed, since $\f$ is generically chaotic, i.e., $C_2(\f)$ is residual in $\mathcal K(I) \times \mathcal K(I)$, hence $\limsup_{n\to \infty} {\rm diam\,} f^n(J)>0$ for any interval $J,$ and thus also 
$$\limsup_{n\to \infty} {\rm diam\,}  f^n({\rm cl}(J)) > 0.$$ 
By Lemma \ref{lm7}, ${\rm Orb}(f, {\rm cl}(J))$ contains a periodic point of $f$ with period $1$ or $2$, hence ${\rm Orb}(g, {\rm cl}(J))$ contains a fixed point of $g$, and consequently $g^s({\rm cl}(J))$ contains a fixed point of $g$ for some $s$;  it suffices to take into account that
$$
\limsup_{n \to \infty} {\rm diam\,} g^n(J)=\limsup_{n \to \infty} {\rm diam\,} g^n({\rm cl}(J))=\limsup_{n \to \infty} {\rm diam\,} g^{sn}({\rm cl}(J)).
$$

Let $\mathcal J$ be the collection of all intervals containing a fixed point of $g$. For any fixed point $p$ of $g$, let $\mathcal J(p)$ be a collection of all intervals $J \in \mathcal J$ containing $p$. To prove that $g$ satisfies (f-2) we show that 
$$
\inf \{ \limsup_{n \to \infty} {\rm diam\,} g^n(J): J \in \mathcal J\}>0.
$$

By a contradiction, assume that for every $i \in \mathbb N$ there exists a fixed point $x_i$ of $g$, and $K_i \in \mathcal J(x_i)$ such that 
$$\lim_{i\to \infty} \limsup_{n \to \infty} {\rm diam\,} g^n(K_i)=0.$$

Without loss of generality, suppose that the sequence of fixed points $x_i$ converges  to a fixed point $p$ of $g$. For any $i\in \mathbb N$, denote 
$$
\limsup_{n \to \infty} {\rm diam\,} g^n(K_i)=\varepsilon_i.
$$

Obviously, for every $i \in \mathbb N$ there is a $k_i \in \mathbb N$ such that for any $k \ge k_i$, we have ${\rm diam\,} g^k(K_i) \le 2\varepsilon_i$. Let $J_i=\bigcap_{k=1}^\infty g^k(K_i)$. Since $x_i \in K_i$ and $g(x_i)=x_i$ we have $x_i \in J_i$ and by the definition of $k_i$ it is obvious that ${\rm diam\  J_i} \leq 4\varepsilon_i$. Moreover, the set $J_i$ is an invariant interval. 

Since the function $\f$ is generically chaotic, by Lemma \ref{new1}(v), the set $C_2(\g)$ is residual. Further, by Lemma \ref{lm12}, $g^k(K_i)$ is not a singleton for every $i$ and $k$.

So we have shown that arbitrarily close to the fixed point $p$ of $g$  there are arbitrarily small $g$-invariant intervals $J_i$. Moreover, by Lemma \ref{new3}, Lemma \ref{lm3}(iii) and by \cite{S1}, Lemma 4.1(i) (i.e., the original version of  our Lemma \ref{new1}(i) for $f$),  for every two intervals $J, J^\prime$, 
$$
\liminf \limits_{n \to \infty} {\rm dist} (g^n(J), g^n(J^\prime))=0.
$$ 
Now it follows from Lemma \ref{lm9} that there exist arbitrarily small $g$-invariant intervals containing the point $p$.  Since, by 
Lemma \ref{new1}(vi), 
$\g$ is generically chaotic, we have a contradiction with Lemma \ref{lm8}. 
\cbd \\

\noindent {\bf Proof of Theorem \ref{main1}.} It follows by Theorem \ref{th1}, Lemmas \ref{lm3} and \ref{new3},  and Lemma \ref{lm16}.
\cbd \\


\section{Densely $\varepsilon$-chaotic $\f$ implies densely $\varepsilon$-chaotic $f$}

Let $(X,d)$ be a compact metric space and $f\in C(X)$.
Recall that a pair $(x,y)\in X \times X$ is called {\it distal} if 
$\liminf_{n\to\infty} d (f^n(x), f^n(y)) > 0$, and {\it asymptotic} if $\lim_{n\to\infty} d (f^n(x), f^n(y)) = 0$; 
it is $\delta$-{\it asymptotic}, with $\delta >0$, if $\limsup_{n\to\infty} d (f^n(x), f^n(y)) <\delta$. 

And for completeness recall also a statement of the well known\\

\noindent{\bf Baire Category Theorem. }{\em
If a non-empty complete metric space is the union of a sequence of its closed subsets, then at least one set in the sequence must have non-empty interior.
}

\begin{lemma}\label{distal}
Let $f\in C(X)$. If the induced map $\f$ is densely chaotic then the set of distal pairs of $f$ is a first category subset of $X\times X$.
\end{lemma}

\proof
Denote by $D\subset X\times X$ the set of distal pairs for $f$. For every $\delta > 0$ and every $n\in \mathbb N$ let
$$
D_{\delta, n} := \{ (x,y)\in X\times X;\ d(f^j(x), f^j (y)) \geq \delta,\ {\rm for\ any\ } j\geq n \}.
$$
Obviously, every $D_{\delta, n}$ is a closed set and $D =\bigcup_{\delta > 0} \bigcup_{n\in\mathbb N} D_{\delta, n} = \bigcup_{k\in\mathbb N} \bigcup_{n\in\mathbb N} \\ D_{1/k, n}$. So, it suffices to show that every $D_{\delta, n}$ is nowhere dense in $X\times X$. Assume on the contrary that, for some $\delta_0 > 0$ and $n_0\in\mathbb N$,
int($D_{\delta_0, n_0})\neq \emptyset$.
Then there are non-empty open sets $U, V \subseteq X$ such that $U \times V \subseteq D_{\delta_0, n_0}$. 
It follows that for every  $j\geq n_0$ and every pair $(x,y) \in f^j(U\times V)$, $d (x,y) \geq \delta_0$,
whence, for every non-empty compact sets $M\subseteq U$, $N\subseteq V$, 
$d_H(\f^j(M), \f^j(N)) \geq \delta_0$ -- a contradiction. 
\cbd \\

\begin{theorem}\label{main2}
Let $f \in C(X)$ be such that the induced map $\f$ is densely $\varepsilon$-chaotic. Then $f$ is densely $\varepsilon$-chaotic.
\end{theorem}

\proof
Assume on the contrary that $f$ is not densely $\varepsilon$-chaotic. Then there are open sets $\emptyset\ne U_0, V_0\subseteq X$ such that  $U_0\times V_0$ contains no $\varepsilon$-Li-Yorke pair. 
Let $D$ the set of distal pairs of $f$, and let for every $n\in\mathbb N$, 
$$
A_{n}:=\{ (x,y)\in U_0\times V_0;\  d(f^j(x),f^j(y))\le\varepsilon, \ {\rm for\ any} \ j\ge n\}.
$$

Obviously, every $A_{n}$ is a closed set, $A_{1}\subseteq A_{2}\subseteq\dots$, and  the set of $\varepsilon$-asymptotic pairs of $f$ contained in  $U_0\times V_0$ is a subset of
$\bigcup_{n\in\mathbb N} A_{n}$. 
Since $U_0\times V_0\subseteq  D \cup \bigcup_{n\in\mathbb N} A_{n}$,
and $D$ is (by Lemma \ref{distal}) a first category set, the Baire Category Theorem implies that  there is an $n_0\in\mathbb N$ such that $A_{n_0}$ has non-empty interior. Consequently, there are non-empty open sets $U_1\subseteq U_0, V_1\subseteq V_0$ such that $U_1\times V_1\subseteq A_{n_0}$. Since $\f$ is densely $\varepsilon$-chaotic there are compact sets $M\subset U_1, N\subset V_1$ such that $\limsup_{j\to\infty}$ $d_H(\overline f^j(M), \overline f^j(N))>\varepsilon$. On the other hand, 
$$
d(f^j(x),f^j(y))\le \varepsilon, \ {\rm for\ every}\ (x,y)\in M\times N \  {\rm and\ every} \  j\ge n_0,
$$ 
whence $\limsup_{j\to\infty}d_H(\f^j(M),\f^j(N))\le\varepsilon$ -- a contradiction.
\cbd 

\bigskip

Using Theorem \ref{main2} and \cite{Mur} we obtain the following
\begin{corollary}\label{cor3}
Let $f \in C(X)$. If the induced map $\f$ is generically $\varepsilon$-chaotic, then $f$ is generically $\varepsilon$-chaotic.
\end{corollary}


\section{Densely ($\varepsilon$-)chaotic $f$ implies densely ($\varepsilon$-)chaotic $\f$} 

In this section we show firstly, that densely chaotic $f$ implies densely chaotic $\f$
(Theorem \ref{main3a}), and secondly, that analogous result is true for  dense $\varepsilon$-chaos
 (Theorem \ref{main3b}).

\begin{theorem}\label{main3a}
Let $f \in C(I)$ be densely chaotic and let $\f$ be induced by $f$. Then $\f$ is densely chaotic.
\end{theorem}

\proof 
Let $\langle U_1, U_2, \ldots , U_n\rangle$ and $\langle V_1, V_2, \ldots , V_m\rangle$ be arbitrary open sets in $\mathcal K(I)$ with $m, n \in \mathbb N$. To prove that $\f$ is densely chaotic it suffices to show that  $\langle U_1, U_2, \ldots , U_n\rangle \times \langle V_1, V_2, \ldots , V_m\rangle$ contains a Li-Yorke pair $(U, V)$ for $\f$ (i.e., that $(U, V) \in C(\f)$). Since we will frequently work  with the systems $\{U_i\, ; i=1, \ldots , n\}$ and $\{V_j\, ; j=1, \ldots , m\}$ in the following proof, we will use, to shorten the notation, the symbols $i$ and $j$ {\it only} in this sense, i.e., $i=1, \ldots, n$ and $j=1, \ldots, m$.
Let $\delta_1:=\min\limits_i \{\liminf\limits_{k\to\infty} {\rm diam\,}  f^k(U_i)\}$ and $\delta_2:=\min\limits_j\{\liminf\limits_{k\to\infty} {\rm diam\,}  f^k(V_j)\}$. By Theorem \ref{th2}, $\delta_1,\delta_2 > 0$.
Put
$$
\delta := \min \{\delta_1,\delta_2 \}.
$$

By Theorem \ref{th2} (a), (b) there are fixed point $x_0 \in I$ of $f$ and  $k \in \mathbb N$ such that for every $i, j$ we have
\begin{equation}\label{eq11}
{\rm dist} (f^k(U_i), x_0)< \delta/4 \quad  {\rm and\quad } {\rm dist} (f^k(V_j), x_0)< \delta/4 ,
\end{equation}
and also
\begin{equation}\label{eq111}
{\rm diam\,} f^k(U_i) > \delta/2 \quad  {\rm and \quad } {\rm diam\,} f^k(V_j) > \delta/2. 
\end{equation}
Denote
$$
{\mathscr S_1}:=\{U_i: x_0-\delta/4 \in f^k(U_i)\} \cup \{V_j: x_0 - \delta/4 \in f^k(V_j)\},
$$
$$
{\mathscr S_2}:=\{U_i: U_i \notin {\mathscr S_1}\} \cup \{V_j: V_j \notin {\mathscr S_1}\}
$$
and put
\begin{equation}\label{eq1}
S_1: = \bigcap_{A \in {\mathscr S_1}}f^k(A),\quad  S_2 := \bigcap_{A \in {\mathscr S_2}}f^k(A).
\end{equation}
Note that, for any set $A \in {\mathscr S_2}$, $x_0+\delta/4 \in f^k(A)$.  

Since a set $A\in {\mathscr S_1}$ is an open interval, by (\ref{eq11}) and (\ref{eq111}) we get that $f^k(A)$ is a non-degenerate interval containing $x_0 - \delta/4$  together with its right neighbourhood and hence, $S_1$ is also a non-degenerate interval. Analogously, $S_2$ is a non-degenerate interval. There are three possibilities.\\

\noindent {\bf Case 1.} There exist $i_1, i_2 \in \{1, \ldots , n\}$ and $j_1, j_2\in \{1, \ldots , m\}$ such that $U_{i_1}, V_{j_1}\in {\mathscr S_1}$ and $U_{i_2}, V_{j_2}\in {\mathscr S_2}$. Then $S_1,$ $S_2$ (see (\ref{eq1})) are non-empty and, since $f$ is densely chaotic, there is a Li-Yorke pair $(s, p) \in S_1 \times S_1$. Let $r \in S_2$ be arbitrary.  For any $i, j$, define $u_i \in U_i$ and  $v_j \in V_j$ in the following way

$$
u_i= \left\{
\begin{array}{l}
f^{-k}(s) \quad {\rm if}\ U_i \in {\mathscr S_1}, \\
f^{-k}(r) \quad {\rm if}\ U_i \in {\mathscr S_2},
\end{array}
\right.
$$

$$
v_j= \left\{
\begin{array}{l}
f^{-k}(p) \quad {\rm if}\ V_j \in {\mathscr S_1}, \\
f^{-k}(r) \quad {\rm if}\ V_j \in {\mathscr S_2}.
\end{array}
\right.
$$
Let $U=\{u_1, u_2, \ldots , u_n\}$ and $V=\{v_1, v_2, \ldots , v_m\}$. Consequently, $U \in \langle U_1, U_2, \ldots , U_n \rangle$, $V \in \langle V_1, V_2, \ldots , V_m \rangle$, and moreover $\f^k(U)=\{s, r\}$ and $\f^k(V) =\{p, r\}$. Obviously, $(U, V)$ is a Li-Yorke pair for $\f$.\\

\noindent {\bf Case 2.} Let ${\mathscr S_1}, {\mathscr S_2}$ be non-empty and assume that one of them contains either all $U_i$, $i=1, \ldots ,n$, or all $V_j$, $j=1, \ldots ,m$. Without loss of generality, assume that $U_i \in {\mathscr S_1}$ for any $i$ (hence ${\mathscr S_2}$ contains only the sets $V_j$ for some $j$; note that ${\mathscr S_1}$ can also contain some intervals $V_j$). Similarly as in Case 1, there exists a Li-Yorke pair $(s, r) \in S_1 \times S_2$. For any $i, j$, define $u_i \in U_i$ and  $v_j \in V_j$ in the following way

$$
u_i= f^{-k}(s)
$$
and
$$
v_j= \left\{
\begin{array}{l}
f^{-k}(s) \quad {\rm if}\ V_j \in {\mathscr S_1}, \\
f^{-k}(r) \quad {\rm if}\ V_j \in {\mathscr S_2}.
\end{array}
\right.
$$
Let $U=\{u_1, u_2, \ldots , u_n\}$ and $V=\{v_1, v_2, \ldots , v_m\}$. Then $\f^k(U)=\{s\}$ and $\{r\} \subseteq \f^k(V) \subseteq \{s, r\}$ and hence $(U, V)$ is a Li-Yorke pair for $\f$.\\

\noindent {\bf  Case 3.} Either ${\mathscr S_1}$ or ${\mathscr S_2}$ is empty. Without loss of generality assume that ${\mathscr S_2} = \emptyset$. Since $f$ is densely chaotic, there is a Li-Yorke pair $(s, p) \in S_1 \times S_1$. For any $i, j$, define $u_i \in U_i$ and  $v_j \in V_j$ by 

$$
u_i= f^{-k}(s) \qquad {\rm and}\qquad v_j=f^{-k}(p).
$$
Let $U=\{u_1, u_2, \ldots , u_n\}$ and $V=\{v_1, v_2, \ldots , v_m\}$. Then $\f^k(U)=\{s\}$ and $\f^k(V) = \{p\}$ and so, $(U, V)$ is a Li-Yorke pair for $\f$.
\cbd \\

The following two results on transitive maps will be used in the proof of Theorem \ref{main3b}.

\begin{lemma} \label{bc} {\rm (\cite{BC}, p. 156, Proposition 42)} 
Let $f: I \rightarrow I$ be transitive. Then exactly one of the following alternatives holds:
\begin{enumerate}
\item[(i)] for every positive integer $s$, $f^s$ is transitive,
\item[(ii)] there exist non-degenerate closed intervals $J, K$ with $J \cup K = I$ and $J \cap K = \{y\}$, where $y$ is a fixed point of $f$, such that $f(J)=K$ and $f(K)=J$. 
\end{enumerate}
\end{lemma}

\begin{lemma}\label{bc2}{\rm (\cite{BC}, p. 157, Proposition 44)} Let $f\in C(I)$. Then $f^2$ is transitive if and only if, for every open subinterval $J$ and every closed subinterval $H$ which does not contain an endpoint of $I$, there is a positive integer $N$ such that $H \subseteq f^n(J)$ for every $n>N$. 
\end{lemma}

\begin{theorem}\label{main3b}
Let $f \in C(I)$ be densely $\varepsilon$-chaotic and let $\f$ be induced by $f$. Then $\f$ is densely $\varepsilon$-chaotic.
\end{theorem} 

\proof
By Theorem \ref{th1}, the map $f$ is generically $\varepsilon$-chaotic. Let $\langle U_1, U_2, \ldots, \\ U_n \rangle$ and $\langle V_1, V_2, \ldots , V_m \rangle$ be open sets in $\mathcal K(I)$, where $n, m \in \mathbb N$ and $U_i, V_j$ are open intervals in $I$ for $1 \le i \le n, 1 \le j \le m$. To prove that $\f$ is densely $\varepsilon$-chaotic it suffices to find an $\varepsilon$-Li-Yorke pair $(U, V) \in \langle U_1, U_2, \ldots , U_n\rangle \times \langle V_1, V_2, \ldots , V_m\rangle$ (i.e., $(U, V) \in C(\f, \varepsilon)$). Similarly as in the proof of Theorem  \ref{main3a}, we will use the symbols $i$ and $j$ {\it only} for this purpose, i.e., $i=1, \ldots, n$ and $j=1, \ldots, m$.

By Theorem \ref{th1} (h-2) and (g-2), there are a positive integer $l$ and open intervals $U'_i, V'_j$, for any $i, j$, such that  $U'_i \subset f^l(U_i) \cap\ {\rm int}(T_{U_i})$ (resp. $V'_j \subset f^l(V_j) \cap\ {\rm int}(T_{V_j})$ ), where $T_{U_i}$ (resp. $T_{V_j}$) is an invariant transitive interval. 

By Theorem \ref{th1} (h-1), $f$ has a unique invariant transitive interval or two transitive intervals having one point in common. \\

\noindent {\bf (a)} Consider first the existence of two invariant transitive intervals $T_1, T_2$  with one common point $x_0$. Obviously, $x_0$ is a fixed point. Assume that $T_1$ lies on the left of $T_2$.

Denote
$$
{\mathscr S_1}:= \{U'_i: U'_i \subset T_1\} \cup \{V'_j: V'_j \subset T_1\},
$$
$$
{\mathscr S_2}:= \{U'_i: U'_i \subset T_2\} \cup \{V'_j: V'_j \subset T_2\}.
$$
By Theorem \ref{th1} (g-2), there is a positive $b$ such that
$\liminf\limits_{n \to \infty} {\rm diam\,}f^n(J)>b$ for every interval $J$, and by Theorem \ref{th1} (g-1), there is a non-negative integer $k$ such that, for every $i, j$, we have
\begin{equation}\label{eq22}
{\rm dist} (f^k(U'_i), x_0)< b/4 \quad  {\rm and\quad } {\rm dist} (f^k(V'_j), x_0)< b/4 ,
\end{equation}
and moreover
\begin{equation}\label{eq222}
{\rm diam\,} f^k(U'_i) > b/2 \quad  {\rm and \quad } {\rm diam\,} f^k(V'_j) > b/2. 
\end{equation}
Put
\begin{equation}\label{eq2}
S_1: = \bigcap_{A \in {\mathscr S_1}}f^k(A),\quad  S_2 := \bigcap_{A \in {\mathscr S_2}}f^k(A).
\end{equation}
Since a set $A\in {\mathscr S_1}$ is an open interval, by (\ref{eq22}) and (\ref{eq222}), $f^k(A)$ is a non-degenerate interval containing $x_0 - b/4$  together with its right neighbourhood, and hence $S_1$ is also a non-degenerate interval. Analogously, $S_2$ is a non-degenerate interval. Similarly as in the proof of Theorem \ref{main3a} there are three possibilities.  \\

\noindent {\bf Case 1.} There exist $i_1, i_2 \in \{1, \ldots , n\}$ and $j_1, j_2\in \{1, \ldots , m\}$ such that $U'_{i_1}, V'_{j_1}\in {\mathscr S_1}$ and $U'_{i_2}, V'_{j_2}\in {\mathscr S_2}$. Then $S_1,$ $S_2$ (see (\ref{eq2})) are non-empty and, since $f$ is densely $\varepsilon$-chaotic, there is an $\varepsilon$-Li-Yorke pair $(s, p) \in S_1 \times S_1$. Let $r \in S_2$ be arbitrary.  For any $i, j$, define $u_i \in U_i$ and  $v_j \in V_j$ in the following way

$$
u_i= \left\{
\begin{array}{l}
f^{-(k+l)}(s) \quad {\rm if}\ U'_i \in {\mathscr S_1}, \\
f^{-(k+l)}(r) \quad {\rm if}\ U'_i \in {\mathscr S_2},
\end{array}
\right.
$$

$$
v_j= \left\{
\begin{array}{l}
f^{-(k+l)}(p) \quad {\rm if}\ V'_j \in {\mathscr S_1}, \\
f^{-(k+l)}(r) \quad {\rm if}\ V'_j \in {\mathscr S_2}.
\end{array}
\right.
$$
Let $U=\{u_1, u_2, \ldots , u_n\}$ and $V=\{v_1, v_2, \ldots , v_m\}$. Consequently, $U \in \langle U_1, U_2, \ldots , U_n \rangle$, $V \in \langle V_1, V_2, \ldots , V_m \rangle$ and moreover, $\f^{k+l}(U)=\{s, r\}$ and $\f^{k+l}(V) =\{p, r\}$. Since $r \in T_2$ and $(s, p) \in T_1 \times T_1$ is an $\varepsilon$-Li-Yorke pair, we have  that $(U, V)$ is an $\varepsilon$-Li-Yorke pair for $\f$.\\

\noindent {\bf Case 2.} Let ${\mathscr S_1}, {\mathscr S_2}$ be non-empty and assume that one of them contains either all $U'_i$, $i=1, \ldots ,n$, or all $V'_j$, $j=1, \ldots ,m$. Without loss of generality we may assume that $U'_i \in {\mathscr S_1}$ for any $i$ (hence ${\mathscr S_2}$ contains only the sets $V'_j$ for some $j$; note that ${\mathscr S_1}$ can also contain some intervals $V'_j$). Similarly as in Case 1, there exists an $\varepsilon$-Li-Yorke pair $(s, r) \in S_1 \times S_2$. For any $i, j$, define $u_i \in U_i$ and  $v_j \in V_j$ in the following way

$$
u_i= f^{-(k+l)}(s),
$$
and
$$
v_j= \left\{
\begin{array}{l}
f^{-(k+l)}(s) \quad {\rm if}\ V'_j \in {\mathscr S_1}, \\
f^{-(k+l)}(r) \quad {\rm if}\ V'_j \in {\mathscr S_2}.
\end{array}
\right.
$$
Let $U=\{u_1, u_2, \ldots , u_n\}$ and $V=\{v_1, v_2, \ldots , v_m\}$. Then $\f^{k+l}(U)=\{s\}$ and $\{r\} \subseteq \f^{k+l}(V) \subseteq \{s, r\}$, and hence $(U, V)$ is an $\varepsilon$-Li-Yorke pair for $\f$.\\

\noindent {\bf  Case 3.} Either ${\mathscr S_1}$ or ${\mathscr S_2}$ is empty. Without loss of generality we assume that ${\mathscr S_2} = \emptyset$. Since $f$ is densely $\varepsilon$-chaotic there is an $\varepsilon$-Li-Yorke pair $(s, p) \in S_1 \times S_1$. For any $i, j$, define $u_i \in U_i$ and  $v_j \in V_j$ by 

$$
u_i= f^{-(k+l)}(s) \qquad {\rm and}\qquad v_j=f^{-(k+l)}(p).
$$
Let $U=\{u_1, u_2, \ldots , u_n\}$ and $V=\{v_1, v_2, \ldots , v_m\}$. Then $\f^{k+l}(U)=\{s\}$ and $\f^{k+l}(V) = \{p\}$,
and so $(U, V)$ is an $\varepsilon$-Li-Yorke pair for $\f$.\\

\noindent {\bf (b)} Now we consider that $f$ has a unique invariant transitive interval $T$. By Lemma \ref{bc}, either $f^s|_T$ is transitive for every positive integer $s$, or there exist non-degenerate closed intervals $J, K$ with $J \cup K = T$ and $J \cap K = \{x_0\}$, where $x_0$ is a fixed point of $f|_T$, such that $f|_T(J)=K$ and $f|_T(K)=J$. Let the second  possibility hold, put $g:=f^2$. Note that since $f$ is densely $\varepsilon$-chaotic (and hence also generically chaotic), $g$ is also densely $\varepsilon$-chaotic (see \cite{S1}, Lemma 4.2). By Lemma \ref{new1}(i)  and since 
$$
\limsup\limits_{n \to \infty} {\rm dist} (\g^n(A), \g^n(B))>\varepsilon \Rightarrow \limsup\limits_{n \to \infty} {\rm dist} (\f^n(A), \f^n(B))>\varepsilon,
$$ 
it suffices to prove  that $\g$ is densely $\varepsilon$-chaotic. Since $J, K$ are two $g$-invariant transitive intervals with one common point $x_0$,  we may proceed analogously as in part  (a) of this proof. Hence $\g$ is densely $\varepsilon$-chaotic.

Finally, let $f^s|_T$ be transitive for every positive integer $s$ (specially $f^2|_T$ is transitive), let $U'$ (respectively $V'$) denote the closure of the  convex hull of $\bigcup_i U'_i$ (respectively $\bigcup_j V'_j$). From the construction of $U'_i$ and $V'_j$ we may assume that $U'$ and $V'$ do not contain an endpoint of $T$.  Then, by Lemma \ref{bc2}, there exists a positive integer $k$ such that $f^k(U'_i) \supseteq U'$ and $f^k(V'_j) \supseteq V'$, for every $i, j$. Since $f$ is densely $\varepsilon$-chaotic, there is an $\varepsilon$-Li-Yorke pair $(s, p) \in U' \times V'$. For any $i, j$, define $u_i \in U_i$ and $v_j \in V_j$ in the following way
$$
u_i=f^{-(k+l)}(s) \quad {\rm and}\quad v_j=f^{-(k+l)}(p).
$$
Let $U=\{u_1, u_2, \ldots , u_n\}$ and $V=\{v_1, v_2, \ldots , v_m\}$. Then $\f^{k+l}(U)=\{s\}$ and $\f^{k+l}(V) =\{p\}$. Hence $(U, V)$ is an $\varepsilon$-Li-Yorke pair for $\f$. The proof of the  theorem is complete. 
\cbd

\medskip

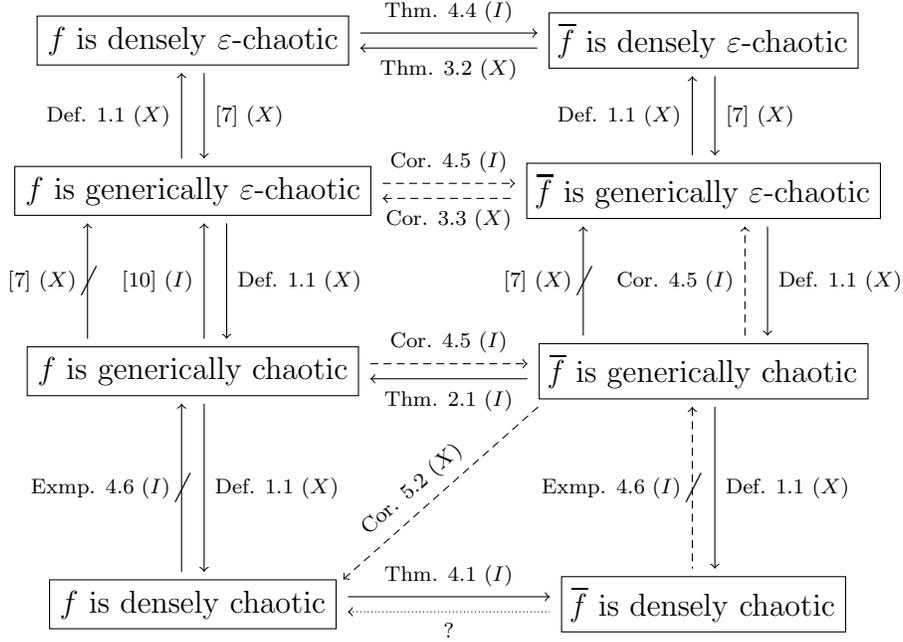
\begin{figure}[h]

\begin{tikzpicture}

\node[draw,rectangle] (uz3) at (-0.09,0.2){$f$ is densely $\varepsilon$-chaotic};
\node[draw,rectangle] (uz33) at (6.7, 0.2) {$\f$ is densely $\varepsilon$-chaotic};
\node[draw,rectangle] (uz2) at (-0.09,-1.8){$f$ is generically $\varepsilon$-chaotic};
\node[draw,rectangle] (uz22) at (6.7,-1.8) {$\f$ is generically $\varepsilon$-chaotic};
\node[draw,rectangle] (uz1) at (-0.09,-4.2){$f$ is generically chaotic};
\node[draw,rectangle] (uz11) at (6.7,-4.2) {$\f$ is generically chaotic};
\node[draw,rectangle] (uz4) at (-0.09,-7.3){$f$ is densely chaotic};
\node[draw,rectangle] (uz44) at (6.7,-7.3) {$\f$ is densely chaotic};



\draw[densely dashed,->] ([xshift=0.15cm, yshift=0.1cm]uz1.east) -- ([xshift=-0.15cm, yshift=0.1cm]uz11.west) node[midway,above] {{\tiny Cor. \ref{cor4} ($I$)}} ;
\draw[<-]  ([xshift=0.15cm, yshift=-0.1cm]uz1.east) -- ([xshift=-0.15cm, yshift=-0.1cm]uz11.west) node[midway,below] {{\tiny Thm. \ref{main1} ($I$)}} ;

\draw[densely dashed, ->] ([xshift=0.15cm, yshift=0.1cm]uz2.east) -- ([xshift=-0.15cm, yshift=0.1cm]uz22.west) node[midway,above] {{\tiny Cor. \ref{cor4} ($I$)}} ; 
\draw[densely dashed, <-] ([xshift=0.15cm, yshift=-0.1cm]uz2.east) -- ([xshift=-0.15cm, yshift=-0.1cm]uz22.west) node[midway,below] {{\tiny Cor. \ref{cor3} ($X$)}} ;

\draw[->] ([xshift=0.15cm, yshift=0.1cm]uz3.east) -- ([xshift=-0.15cm, yshift=0.1cm]uz33.west) node[midway,above] {{\tiny Thm. \ref{main3b} ($I$)}} ;
\draw[<-] ([xshift=0.15cm, yshift=-0.1cm]uz3.east) -- ([xshift=-0.15cm, yshift=-0.1cm]uz33.west) node[midway,below] {{\tiny Thm. \ref{main2} ($X$)}} ;

\draw[->] ([xshift=0.15cm, yshift=0.1cm]uz4.east) -- ([xshift=-0.15cm, yshift=0.1cm]uz44.west) node[midway,above] {{\tiny Thm. \ref{main3a} ($I$)}} ;

\draw[densely dotted,<-] ([xshift=0.15cm, yshift=-0.1cm]uz4.east) -- ([xshift=-0.15cm, yshift=-0.1cm]uz44.west) node[midway,below] {{\tiny ?}} ;


\draw[->] ([xshift=0.15cm, yshift=0.1cm]uz1.north) --  ([xshift=0.15cm, yshift=-0.1cm]uz2.south) node[midway,left] { {\tiny \cite{S1} ($I$)}}; 
\draw[<-] ([xshift=0.45cm, yshift=0.1cm]uz1.north) -- ([xshift=0.45cm, yshift=-0.1cm]uz2.south) node[midway,right] { {\tiny Def. \ref{def} ($X$)}}; 
\draw[->] ([xshift=-1.4cm, yshift=0.1cm]uz1.north) -- node[midway, left] {\(\not\)} ([xshift=-1.4cm, yshift=-0.1cm]uz2.south) node[midway,left] { {\tiny \cite{Mur}  ($X$)}}; 

\draw[->] ([xshift=-0.15cm, yshift=0.1cm]uz2.north) --  ([xshift=-0.15cm, yshift=-0.1cm]uz3.south) node[midway,left] { {\tiny Def. \ref{def} ($X$)}}; 
\draw[<-] ([xshift=0.15cm, yshift=0.1cm]uz2.north) -- ([xshift=0.15cm, yshift=-0.1cm]uz3.south) node[midway,right] { {\tiny \cite{Mur} ($X$)}}; 

\draw[<-] ([xshift=-0.15cm, yshift=-0.1cm]uz1.south) -- node[midway, left] {\(\not\)} ([xshift=-0.15cm, yshift=0.1cm]uz4.north) node[midway,left] {{\tiny Exmp. \ref{exmp} ($I$)}} ;
\draw[->] ([xshift=0.15cm, yshift=-0.1cm]uz1.south) -- ([xshift=0.15cm, yshift=0.1cm]uz4.north)node[midway,right] {{\tiny Def. \ref{def} ($X$)}}; 


\draw[densely dashed, ->] ([xshift=0.55cm, yshift=0.1cm]uz11.north) --  ([xshift=0.55cm, yshift=-0.1cm]uz22.south)node[midway,left] {{\tiny Cor. \ref{cor4} ($I$)}} ;
\draw[<-] ([xshift=0.85cm, yshift=0.1cm]uz11.north) -- ([xshift=0.85cm, yshift=-0.1cm]uz22.south) node[midway,right] {{\tiny Def. \ref{def}  ($X$)}}; 
\draw[->] ([xshift=-1.6cm, yshift=0.1cm]uz11.north) -- node[midway, left] {\(\not\)} ([xshift=-1.6cm, yshift=-0.1cm]uz22.south) node[midway,left] { {\tiny \cite{Mur}  ($X$)}}; 

\draw[->] ([xshift=-0.15cm, yshift=0.1cm]uz22.north) --  ([xshift=-0.15cm, yshift=-0.1cm]uz33.south) node[midway,left] { {\tiny Def. \ref{def} ($X$)}}; 
\draw[<-] ([xshift=0.15cm, yshift=0.1cm]uz22.north) -- ([xshift=0.15cm, yshift=-0.1cm]uz33.south) node[midway,right] { {\tiny \cite{Mur} ($X$)}}; 

\draw[->] ([xshift=0.15cm, yshift=-0.1cm]uz11.south) -- ([xshift=0.15cm, yshift=0.1cm]uz44.north)node[midway,right] {{\tiny Def. \ref{def} ($X$) }}; 
\draw[densely dashed, <-] ([xshift=-0.15cm, yshift=-0.1cm]uz11.south) -- node[midway, left] {\(\not\)} ([xshift=-0.15cm, yshift=0.1cm]uz44.north) node[midway,left] {\tiny Exmp. \ref{exmp} ($I$)} ;

\draw[densely dashed, ->] ([xshift=-2.2cm, yshift=-0.1cm]uz11.south) -- ([xshift=2cm]uz4.north)node[sloped,above, midway]{{\tiny\  Cor. \ref{cor5} ($X$)\ \ \ \ \ \ \ }};

\end{tikzpicture}
\caption{The scheme.}\label{scheme}
\end{figure}

From Theorem \ref{main3b}, \cite{Mur}, \cite{S1} and Theorem \ref{main1} (see also the scheme on Figure \ref{scheme}) 
we obtain the following

\begin{corollary}\label{cor4}
Let $f \in C(I)$ be generically $\varepsilon$-chaotic (resp., generically chaotic), then the induced map $\f$ is generically $\varepsilon$-chaotic (resp., generically chaotic).
Moreover, if $\f$ is generically chaotic, then $\f$ is generically $\varepsilon$-chaotic.
\end{corollary}

\noindent{\bf Remark.}
Note that the equivalence of generic and generic $\varepsilon$-chaos for both $f$ and $\f$ in the interval case is no more true for maps on general compact metric spaces: Murinov\' a in \cite{Mur} gave an example of a map $f$ on the plane which is generically chaotic but not generically $\varepsilon$-chaotic. It is easily seen that the same is true for the induced map $\f$. 

\medskip
Now we will present an example  of a densely chaotic function $f : I \to I$  such that neither $f$, nor $\f$ is densely $\varepsilon$-chaotic, for any $\varepsilon >0$. 
This example is taken from \cite{S1} but, for convenience, we give it here adding its graph.  

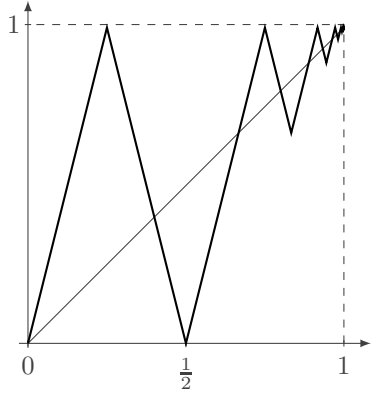
\begin{figure}[h]
\begin{center}
\begin{tikzpicture}[scale=0.6] 

\draw[-latex,color=darkgray,thin] (-0.2,0) -- (7.6,0);
\draw[shift={(0,0)},color=darkgray,thin] (0pt,2pt) -- (0pt,-2pt)
                                   node[below] {\footnotesize $0$};
\draw[shift={(3.5,0)},color=darkgray,thin] (0pt,2pt) -- (0pt,-2pt)
                                   node[below] {\footnotesize $\frac{1}{2}$};
\draw[shift={(7,0)},color=darkgray,thin] (0pt,2pt) -- (0pt,-2pt)
                                   node[below] {\footnotesize $1$};  
                                   
\draw[-latex,color=darkgray,thin] (0,-0.1) -- (0, 7.6);
\draw[shift={(0, 7.069)},color=darkgray,thin] (-2pt,0pt) -- (2pt,0pt)
                                   node[left] {\footnotesize $1$};                                                                      

\draw[dashed,color=darkgray] (7,0)--(7,7.069)--(0, 7.069); 

\draw[thin, color=darkgray] (0,0)--(7,7);  

\draw[thick] (0,0)--(1.75, 7)--(3.5, 0)
	--(4.667, 4.667)--(5.25, 7)--(5.833, 4.667)
	--(6.222, 6.222)--(6.417, 7)--(6.611, 6.222)
	--(6.7407, 6.7407)--(6.8056, 7)--(6.87037, 6.7407)
	--(6.91358, 6.91358)--(6.93519, 7)--(6.95679, 6.91358)
	--(6.9711934, 6.9711934)--(6.978395, 7)--(6.9855967, 6.9711934)
	--(6.990397805, 6.990397805)--(6.99279835390947, 7)--(6.99519890260631, 6.990397805)
	--(6.99679926840421, 6.99679926840421)--(6.99759945130316, 7)--(6.9983996342021, 6.99679926840421)
	--(6.99893308946807, 6.99893308946807)
	--(7,7);

\end{tikzpicture}

\end{center} 

\caption{Sketch of the graph of $f$.} \label{graph}
\end{figure}

\exmp (\cite{S1}, Example 3.6) \label{exmp} For $n = 0, 1, 2, \ldots$, let 
$$
a_n:=1-\frac{1}{3^n}, \quad b_n:= 1-\frac{1}{4\cdot 3^{n-1}},\quad c_n:=1-\frac{1}{2 \cdot 3^{n}}, 
$$
and let $I_n$ be a closed interval $[a_n, 1]$. Define a sequence $\{f_n\}_{n=0}^\infty$ of piecewise linear functions $f_n \in C(I) $ such that for $i=0, 1, \ldots , n$, $f_n$ is linear on each of intervals $[a_i, b_i]$, $[b_i, c_i]$, $[c_i, a_{i+1}]$,
$$
f_n(a_i)=a_i, \quad f_n(b_i)=1, \quad f_n(c_i)=a_i 
$$
and $f_n|_{I_{n+1}}$ is the identity function. Let $f \in C(I)$ be the uniform limit of $f_n$ for $n \to \infty$. A sketch of the graph of $f$ is provided in Figure \ref{graph}. 

The function $f$ is densely chaotic, but not densely $\varepsilon$-chaotic for any $\varepsilon>0$ (for more details, see \cite{S1}). Thus, by Theorems \ref{main3a} and \ref{main1}, $\f$ is densely chaotic, but not densely $\varepsilon$-chaotic for any $\varepsilon>0$.


\section{Related results and an open problem}

The results of this section are motivated by the question of transmission of dense chaos from $\f$ to $f$.
We assume that $(X,d)$ is a compact metric space and $f\in\mathcal C(X)$. We do not know whether dense chaoticity of   $\overline f$ implies the same property for $f$, even for interval maps. Instead, we are able to prove the following weaker result. Its proof is based on similar techniques as that of Theorem \ref{main2}.\\

Let us recall that by  $C(f) = C_1(f)\cap C_2(f) \subseteq X \times X$ 
(resp. $C(\f) = C_1(\f)\cap C_2(\f) \subseteq \mathcal K(X) \times \mathcal K(X)$)
 we denote the set of Li-Yorke pairs for $f$ (resp. $\f$); and similarly for $C(f,\varepsilon)$, resp. $C(\f, \varepsilon)$,
see Section \ref{chaos}.


\begin{theorem}\label{husty}
Let $f\in C(X)$. Let for any non-empty open set $G \subseteq \mathcal K(X) \times \mathcal K(X)$, $G \cap C(\f)$ be a second category set.  Then $f$ has also this property (i.e., for any open $\emptyset \neq G \subseteq X \times X$, $G \cap C(f)$ is a second category set). 

\end{theorem}

\proof
Assume on the contrary that $f$ does not have the required property.  Then there are a first category set $E$,
and open sets  $\emptyset \neq U_0, V_0 \subseteq X$ such that every pair of points in $(U_0 \times V_0) \setminus E$ is  either distal or asymptotic.
Denote by $C_{U\times V}(\f, \varepsilon) := C(\f, \varepsilon) \cap (\mathcal K(U) \times \mathcal K(V))$.
Obviously, $C_{U\times V}(\f, \varepsilon_1) \supset C_{U\times V}(\f, \varepsilon_2)$, whenever  $\varepsilon_1 < \varepsilon_2$.
Since $\bigcup_{k\in\mathbb N} C_{U_0 \times V_0}(\f, \frac 1k)$ is a second category set, there is an $\varepsilon_0 > 0$
and non-empty open sets $U_1 \subseteq U_0$ and  $V_1 \subseteq V_0$ such that 
$C(\f, \varepsilon_0)$  is dense in $\mathcal K(U_1) \times \mathcal K(V_1)$.

Denote by $A$ the set of asymptotic pairs of $f$ contained in $U_1 \times V_1$, and, for every $\delta > 0$ and $n\in\mathbb N$, let
$$ 
A_{\delta, n} = \{ (x,y) \in U_1\times V_1;\ d(f^j(x), f^j(y)) \leq \delta ,\ {\rm for\ any\ } j\geq n \}.
$$
Obviously, for every $\delta > 0$ and $n\in\mathbb N$, $A_{\delta, n}$ is a closed set, 
$A\subseteq \bigcup_{n\in\mathbb N} A_{\delta, n}$, and
$A_{\delta, 1} \subseteq A_{\delta, 2} \subseteq \dots$.
Since $E$ is a first category set, Lemma \ref{distal} and Baire category theorem imply that, for any $\delta > 0$ there is an $n_0 \in \mathbb N$ such that $A_{\delta, n_0}$ has non-empty interior.

Consequently, there are non-empty open sets $U_2 \subseteq U_1$, $V_2 \subseteq V_1$, and $n_0\in\mathbb N$ such that $U_2 \times V_2 \subseteq A_{{\varepsilon_0}/{2}, n_0}$.
Since $C(\f, \varepsilon_0)$  is dense in $\mathcal K(U_2) \times \mathcal K(V_2)$, there are compact sets
$M \subset U_2$, $N\subset V_2$ such that $\limsup_{j\to\infty} d_H(\f^j(M), \f^j(N)) \geq\varepsilon_0$.
On the other hand, 
$$
d(f^j(x), f^j(y)) \leq \frac{\varepsilon_0}{2}, {\rm \ for\ every\ } (x,y) \in M\times N {\rm \ and\ every\ } j\geq n_0,
$$
whence $d_H(\f^j(M), \f^j(N)) < \varepsilon_0 /2$ -- a contradiction.
\cbd

\medskip
Obviously, the property from Theorem \ref{husty} is stronger than dense chaoticity;
 note that the map from Example \ref{exmp} is densely chaotic but does not have this property.
On the other hand, any generically chaotic map has this property, so we get the following

\begin{corollary}\label{cor5}
Let $f\in C(X)$. If the induced map $\f$ is generically chaotic, then $f$ is densely chaotic.
\end{corollary}

Let us now survey the obtained results in a scheme. 
In this scheme the arrows mean implications, the dashed arrows are corollaries which follow by the transitivity of implications.
Since some of the properties are true not only for interval maps but even for maps on general compact metric spaces, this fact together with the numbers of the corresponding theorems and corollaries are indicated in the brackets next to the corresponding arrows.

\medskip
\noindent{\bf Remark.}
Note that for piecewise monotone interval maps $f$ (with finite number of pieces of monotonicity) all four variants of generic or dense chaos we deal with in the present paper are equivalent (see \cite{S2}). 	
It can be easily shown that analogous result is true for $\f$. Really, if $\f$ is densely chaotic, then obviously $f$ fulfills property (b) from Theorem \ref{th2} and (f-1) from Theorem \ref{th1}.
From \cite{S2} it follows that these two properties imply (a) from Theorem \ref{th2}.
Using Theorem 1.4 from \cite{S2} we obtain that $f$ is generically chaotic, and hence also $\f$ is generically chaotic (see Corollary \ref{cor4}).
Conseqently, for piecewise monotone interval maps, generic, generic $\varepsilon$-, dense and 
dense $\varepsilon$-chaos for both $f$ and $\f$ are all equivalent.\\

\noindent{\bf Open problem.} The question whether dense chaos is transmitted from $\f$ to $f$ for  maps on general  compact metric spaces as well as for (non piecewise monotone) interval maps (i.e., the dotted arrow indicated with the question mark in the scheme) remains open.\\

Let us conclude with the proof of the property that dense chaoticity of $\f$ cannot occur in spaces with isolated points.

\begin{lemma}
Let $f\in C(X)$. If the induced map $\f$ is densely chaotic, then the set of asymptotic pairs of $f$ has empty interior.
\end{lemma}

\proof
Let $A$ be the set of asymptotic pairs of $f$ and denote by
$A_{k, n} = \{ (x,y) \in X\times X;\ d(f^j(x), f^j(y)) \leq 1/k ,\ {\rm for\ any\ } j\geq n \}$.
Assume, on a contrary, that $int(A)\neq \emptyset$.

Then there are open sets $\emptyset\neq U, V \subseteq X$ such that $U\times V \subseteq A$ and, by the hypothesis, there are compact sets $M\subseteq U, N\subseteq V$ such that $(M, N) \in C(\f, \varepsilon)$.
Let $k>1/\varepsilon$. Since $\bigcup_{n\in\mathbb N} A_{k,n} \supseteq A$ and, for every $n\in \mathbb N$,
$A_{k,n} \subseteq A_{k, n+1}$, there is an $n_0 \in \mathbb N$ such that $A_{k, n_0} \supseteq M\times N$.
Hence  $d(f^j(x), f^j(y)) < 1/k < \varepsilon$, for every $x\in M$, $y\in N$, and any $j\geq n_0$.
Consequently, $\limsup_{j\to\infty} d_H(\f^j(M), \f^j(N)) < \varepsilon$ -- a contradiction.
\cbd \\

\begin{corollary}
Let $f\in C(X)$. If the induced map $\f$ is densely chaotic, then $X$ has no isolated points.
\end{corollary}

At the very end, we give an example of a system with an interesting behavior - while the system itself is asymptotic, the induced system contains LY-pairs. 
\exmp \label{exmp2} Let $\Sigma_2=\{0, 1\}^\mathbb N$ be the space of all sequences of two symbols $0$ and $1$ equipped with the metric $d(x, y)=\max\{1/i; x_i \neq y_i\}$, for any distinct $x=\{x_i\}_{i\in\mathbb N}$ and $y=\{y_i\}_{i\in \mathbb N}$. Denote by $A$ the set of all sequences containing at most one symbol $0$ (i.e., the sequence $1^\infty$ of ones and all sequences of the form $1^r01^\infty$ where $r$ is a nonnegative integer). Let $\sigma : A \rightarrow A$ be the standard shift map. 

Obviously, $\sigma(A)=A$ and for every $x \in A$ there is $n \in \mathbb N$ such that $\sigma^n(x)= 1^\infty$. Hence any pair $(x, y) \in A \times A$ is asymptotic and the set of LY-pairs $C(\sigma)=\emptyset$. 
On the other hand, for the induced system, $C(\overline{\sigma})\neq\emptyset$. Indeed, if  $M:=\{1^\infty\}$ and
 $$N:=\{1^{n_i}01^\infty ; i=0, 1, 2, \ldots, n_0=0\ {\rm and}\ n_{i+1}=n_i+i+2 \}$$
then $(M, N) \in \mathcal K(A)\times \mathcal K(A)$ is a LY-pair for $\overline{\sigma}$.  

Note that any LY-pair on a shift space is always an $\varepsilon$-LY-pair, for some $\varepsilon>0$ (in our space $(\Sigma_2, d)$, $\varepsilon = 1$).\\

\noindent{\bf Acknowledgments.} The authors would like to thank prof. J. Sm\' \i tal for fruitful discussions and valuable comments. 

\end{document}